\documentclass[a4paper,11pt]{article}
\usepackage{amsfonts, amssymb, amscd, amsmath}

\numberwithin{equation}{section}

\newtheorem{theo}{Theorem}[section]
\newtheorem{coro}[theo]{Corollary}
\newtheorem{lemm}[theo]{Lemma}
\newtheorem{prop}[theo]{Proposition}
\newtheorem{defi}[theo]{Definition}
\newtheorem{rema}[theo]{Remark}
\newtheorem{exam}[theo]{Example}

\newenvironment{proo}{\noindent \textbf{{Proof.}} \sf}

\def\qed{\hfill $\diamond$ \bigskip}

\newcommand{\Bdf}{\begin{defi}}
\newcommand{\Edf}{\end{defi}}
\newcommand{\Bte}{\begin{theo}}
\newcommand{\Ete}{\end{theo}}
\newcommand{\Bpo}{\begin{prop}}
\newcommand{\Epo}{\end{prop}}
\newcommand{\Bcr}{\begin{coro}}
\newcommand{\Ecr}{\end{coro}}
\newcommand{\Blm}{\begin{lemm}}
\newcommand{\Elm}{\end{lemm}}
\newcommand{\Bex}{\begin{exam}}
\newcommand{\Eex}{\end{exam}}
\newcommand{\Bdm}{\begin{proo}}
\newcommand{\Edm}{\end{proo}}
\newcommand{\Brm}{\begin{rema}}
\newcommand{\Erm}{\end{rema}}
\newcommand{\hra}{\hookrightarrow}
\newcommand{\ra}{\rightarrow}

\def\B{{\mathcal B}}
\def\C{{\mathcal C}}
\def\D{{\mathcal D}}
\def\E{{\mathcal E}}

\def\N{{\mathcal N}}

\begin{document}
\sf
\title{\sf Resolving by a free action linear category and applications to Hochschild-Mitchell (co)homology}

\author{Claude Cibils and Eduardo N. Marcos
\thanks{\sf The authors were partially supported by the project USP-Cofecub. The second  mentioned author was supported by the thematic project of FAPESP 2014/09310-5, a research grant from CNPq 302003/2018-5. We acknowledge support from the ``Brazilian-French Network in Mathematics".}}
\date{}
\maketitle

\begin{abstract}

Let $G$ be a group acting on a small category $\C$ over a field $k$, that is $\C$ is a $G$-$k$-category. We first obtain an unexpected result: $\C$ is resolvable by a category which is $G$-$k$-equivalent to it, on which $G$ acts freely on objects.

This resolving category enables to show that if the coinvariants and the invariants functors are exact, then the coinvariants and invariants of the Hochschild-Mitchell (co)homology of $\C$ are isomorphic to the trivial component of the Hochschild-Mitchell (co)ho\-mo\-logy of the skew category $\C[G]$.

If the action of $G$ is free on objects, then there is a canonical decomposition of the Hochschild-Mitchell (co)ho\-mo\-lo\-gy of the quotient category $\C/G$ along the conjugacy classes of $G$. This way we provide a general frame for monomorphisms which have been described previously in low degrees.

\end{abstract}

\noindent 2020 MSC: 16E40,  18D20, 16E45
\renewcommand{\contentsname}{\sf Index}.
\tableofcontents

\section{\sf Introduction}

Let $k$ be a field. A $k$-category $\C$ is a small category enhanced over the category of $k$-vector spaces. B. Mitchell in \cite{mitchell} called these categories ``algebras with several objects". Indeed, a $k$-category $\C$ with only a finite set of objets provides a $k$-algebra $a(\C)$ through the direct sum of all its morphisms. B. Mitchell introduced Hochschild-Mitchell (co)homology of $\C$,  see also \cite{cr,gp,keller, mccarthy}. If the number of objects of $\C$ is finite, then the Hochschild-Mitchell (co)homology of $\C$ is isomorphic to the Hochschild (co)homology of $a(\C)$, see  for instance \cite{cr}.

Let $G$ be a group. A $G$-$k$-category is a $k$-category $\C$ with an action of $G$ on $\C$. More precisely there is a group homomorphism
$G\ra \mathsf{Aut}_k\C$, where $\mathsf{Aut}_k\C$ is the group of $k$-functors $\C \to \C$ which are isomorphisms. In this situation there exists  a  skew category $\C[G]$. If the number of objects of $\C$ is finite, then $a(\C[G])$ is  isomorphic to the usual skew group algebra $a(\C)[G]$, see  \cite{cm}.

The action of $G$ on $\C$ is called free if the action on the objects of $\C$ is free.
In this case the quotient $k$-category $\C/G$ exists,
see for instance \cite{cr}. Moreover the functor $\C \to \C/G$ is a Galois covering. This construction has several applications in representation theory,
see for instance \cite{bg, ri, gabriel,abm,kp}.

A central result obtained in \cite{cm} is that if the action of $G$ is free, then  $\C/G$ and  $\C[G]$ are equivalent $k$-categories. However $\C[G]$ exists for any action, while a free action is essential for defining $C/G$.  Hence for a non free action $\C[G]$ is a substitute of $\C/G$.

In the same vein, in this paper we first introduce in Section \ref{section resolving} a resolving main result. Namely let $\C$ be a $G$-$k$-category. There exists a $G$-$k$-category $M_G(\C)$ - called the resolving category - which has a free action of $G$ and which is  $G$-$k$-equivalent to $\C.$ The resolving category is related to the infinite matrix algebra  considered by J. Cornick in \cite{cornick}, which in turn is linked with Cohen-Montgomery duality in \cite{com}, see also \cite{cs}.

In Sections \ref{HM homology} and \ref{HM cohomology} we will use the resolving category for Hochschild-Mitchell homology and cohomology. Indeed, for a $G$-$k$-category $\C$ the resolving category enables to relate the Hoch\-schild-Mitchell (co)homologies of $\C$, of $\C/G$, and of $\C[G]$. In doing so, we underline that we do not change the
bimodules of coefficients, in contrast with the Cartan-Leray type spectral sequence obtained in \cite{cr}, see also \cite{hsjpaa} and \cite{pr}.

In Subsection \ref{section Gk and graded}  we show that there is a direct sum  decomposition of the Hochschild-Mitchell homology of a $G$-graded $k$-category along the set of conjugacy classes of $G$, as in the case of $G$-graded algebras, see \cite{cornick, lorenz, stefan}.

For instance the homology of the graded category $\C[G]$ decomposes as above. Suppose that the action is free and the coinvariants functor is exact. In Subsection \ref{section homology skew} we prove that there is an isomorphism between the trivial conjugacy class direct summand of the homology of $\C[G]$, and the coinvariants of the homology of $\C$.

The resolving category and results of E. Herscovich in \cite{herscovich} enables then to prove that the isomorphism quoted above also holds when the action is not free.

In Subsection \ref{section Galois}  we infer the analog results for a Galois covering.

For Hochschild-Mitchell cohomology, the defining cochains are direct products of vector spaces while for homology the chains are direct sums. Nevertheless, for a graded $k$-category, we get to show in Subsection \ref{section cohomology graded} that the complex of cochains is also a direct product along the conjugacy classes of $G$. Moreover, the subcomplex associated to the trivial conjugacy class is a subdifferential graded algebra of it.

If the action is free and if the invariants functor is exact, we show that there is an isomorphism of algebras between the trivial component
of the cohomology of $\C[G]$ and the invariants of the cohomology of $\C$. Note that the proof is quite different than the one for homology. The resolving category enables then to extend the result for a non free action.

In Subsection \ref{section Galois cohomology} we translate these results  into the Galois covering setting. An immediate consequence is that the invariants of the cohomology of $\C$ are a canonical direct summand of the cohomology of $C/G$. We recover this way the monomorphism obtained in \cite{mmm} which is made explicit in low
  degrees in \cite{ghs}.

   In Section \ref{section skew group algebras} we restrict our results to the case of a $k$-algebra with a finite group acting by automorphism, and its Hochschild cohomology. We underline that the proof relies in an essential way on the existence of the resolving category.

   When specialising to algebras, the results of this work are related with explicit computations made for instance in \cite{farinati, gk,sw,nsw} for Hochschild
   (co)homology of specific skew group algebras, in particular for the symmetric algebra over a finite dimensional vector
   space $V$ over a field $k$, with $G$ a finite subgroup of $GL(V)$ which order is invertible in $k.$

Finally we observe that the results obtained in this work rely on the standard complexes of chains and of cochains which compute the Hochschild-Mitchell (co)homology of a $k$-category. It is well known that these complexes are generally not suitable for effective computations. A similar work on reduced or normalised resolutions will allow to consider concrete situations where our theoretical results may apply. This could be a subject of another paper.

 \textbf{Acknowledgements:} we warmly thank the referee for the very efficient work done.

\section{\sf The resolving category with free action}\label{section resolving}

Let $k$ be a field. A $k$-category is a small category $\C$ enriched over the category of $k$-vector spaces. In other words the objects of $\C$ are a set denoted $\C_0$, for any pair of objects $x, y\in \C_0$
the space of morphisms $_y\C_x$ from $x$ to $y$ has a structure of vector space, the composition in $\C$  is  $k$-bilinear, and the image of the canonical inclusion $k \hra  {}_x\C_x$ is central in $_x\C_x$ for all $x\in \C_0$. In particular $_x\C_x$ is a
$k$-algebra for any $x \in \C_0$.

We often write  $_yf_x$ for a morphism $f$ from $x$ to $y$, that is belonging to $_y\C_x$.

\Bdf
Let $G$ be a group. A \emph{$G$-$k$-category}  is a $k$-category $\C$ with an action of $G$ by $k$-isomorphisms of $\C$.
\Edf

\begin{rema} Equivalently, a $G$-$k$-category $\C$ is a $k$-category  with an action of $G$ on the set  of objects $\C_0$ such that for all $s\in G$ and $_yf_x \in {}_yC_x$ there is an element $_{sy}(sf)_{sx} \in{}_{sy}\C_{sx}$. Moreover, the  map  $_y\C_x \to _{sy}\!\C_{sx}$ given by $f\mapsto sf$ is $k$-linear. For all $s,t\in G$ and for any morphism $f$ we have $t(sf) =(ts)f$. Finally for any object $x$ and any $s\in G$, we have $s(_x1_x)= {}_{sx} 1_{sx}$.
\end{rema}

\begin{exam}\label{single object}
Let $\Lambda$ be a $k$-algebra and let $\Lambda_1$ be the single object $k$-category where the endomorphism algebra of the object is $\Lambda$. Let $G$ be a group acting by algebra automorphisms of $\Lambda$. Then $\Lambda_1$ is a $G$-$k$-category.
\end{exam}

\Bdf
A $G$-$k$-category $\C$ is  a category \emph{with a free action of $G$} if the action of $G$ on its objects is free. That is if $s\in G$ and $x\in \C_0$, then $sx=x$ only holds if $s=1$.
\Edf

\begin{rema}
Except when $G$ is trivial, there is no free action of $G$ on the $G$-$k$-category $\Lambda_1$ of Example \ref{single object}.
\end{rema}

Recall that for a finite object $k$-category $\C$, its $k$-algebra is $$a(\C) = \displaystyle{\bigoplus_{x,y \in \C_0}\  _y\C_x}$$
with product given by matrix multiplication combined with the composition of $\C$. The identity element of $a(\C)$ is the sum of the identities of the objects.

\Brm\label{notalgebramaps}
Notice that if $\C$ and $\D$ are finite object $k$-categories, a $k$-functor $F:\C \ra \D$  does not provides in general a multiplicative $k$-morphism between $a(\C)$ and $a(\D)$.
\Erm

The following result is a straightforward generalisation of Example \ref{single object}.
\Blm
Let  $\C$  be a $G$-$k$-category with a finite set of objets. The group $G$ acts on $a(\C)$ by algebra automorphisms.
\Elm

Next we will show that a $G$-$k$-category $\C$ is resolvable by a $G$-$k$-category with a free action of $G$. That is the resolving category is $G$-$k$-equivalent to $\C$.

\Bdf \label{MGC} (see also \cite{cornick}, \cite{coell})
Let $\C$ be a $G$-$k$-category.
The objects of the \emph{resolving $k$-category} $M_G(\C)$ are $G\times \C_0$. The $k$-vector space of morphisms of $M_G(\C)$ from $(s,x)$ to $(t,y)$ is
$$_{(t,y)}(M_G(\C))_{(s,x)} = {}_y\C_x.$$
The composition is given in the evident way by the composition of $\C.$ The action of $G$ on  $M_G(\C)$ is defined as follows: for $r\in G$, let $r(s,x) = (rs,rx),$ and for
$$f\in {}_{(t,y)}(M_G(\C))_{(s,x)} = {}_{y}\C_{x} \text{ let } rf \in {}_{ry}\C_{rx}={}_{(rt, ry)}(M_G(\C))_{(rs,rx)}.$$
\Edf
Notice that the above action of $G$ is free on the objects of $M_G(\C)$.

\begin{exam}
 \label{resolvant of an algebra}
Let $\Lambda$ be a $k$-algebra with a group acting by automorphisms, and let $\Lambda_1$ be its single object $G$-$k$-category.

\begin{itemize}
\item The resolving category $M_G(\Lambda_1)$ has set of objects $G$.

\item Each vector space of morphisms of $M_G(\Lambda_1)$ is a copy of $\Lambda$.

\item The composition of $M_G(\Lambda_1)$ is given  by the product of $\Lambda$.

\item The action of $G$ on the objects is the product of $G$. On morphisms it is given by the action on objects followed by the action on $\Lambda$.
    \end{itemize}
\end{exam}

Next we describe  the resolving category of $\C$ as a tensor product of categories.
\Bdf
Let $\C$ and $\D$ be $k$-categories. Their \emph{tensor product} $\C\otimes \D$ has set of objects $\C_0 \times \D_0$. Its morphisms are given by:
$$_{(c', d')}(\C\otimes \D)_{(c,d)}={}_{c'}\\C_c\otimes {}_{d'}\D_{d}$$ with the obvious composition.

\Edf
If $\C$ and $\D$ are $G$-$k$-categories then $\C\otimes \D$ is a $G$-$k$-category through the diagonal action of $G$.

\Brm
Let $\C$ be a $G$-$k$-category. Let $k_1$ be the $G$-$k$-trivial category, with one object whose endomorphisms are given by $k$ and trivial $G$-action. We have $M_G(\C) = M_G(k_1)\otimes\C.$
\Erm

\Bte\label{L}
Let $\C$ be a $G$-$k$-category. There is an equivalence of $G$-$k$-categories $L:M_G(\C) \to \C.$
\Ete
\Bdm
 Let $L: M_G(\C) \to \C$ be the
functor defined on the objects by $L(s,x) = x$, while on morphisms $L$ is given by the suitable identity maps. Hence $L$ is a fully faithful $G$-functor which is surjective on the objects, so it is an equivalence of $G$-$k$-categories. \qed
\Edm

We end this section by describing the algebra associated to a resolving category.

\Bpo\label{resolving algebra of an algebra}
Let $\Lambda$ be a $k$-algebra with an action of a finite group $G$ by automorphisms of $\Lambda$. Let $M_G(\Lambda_1)$ be the resolving category of $\Lambda_1$.

The $k$-algebra $a(M_G(\Lambda_1)$  is isomorphic to the matrix algebra $M_G(\Lambda)$ with columns and rows indexed by $G$. The action of $G$ on a matrix is the combination of the action of $G$ on $\Lambda$ with the action of $G$ on the indices of the rows and of the columns.\Epo

The proof relies on the description given at example \ref{resolvant of an algebra}.

\begin{rema}\label{free idempotents}
Let $E$ be the set of diagonal idempotents of the algebra above $M_G(\Lambda)$, each one is given by $1\in\Lambda$ on a certain spot of the diagonal and zeros elsewhere. The number of elements of $E$ is $\vert G\vert$.

The group $G$ acts freely on $E$.
\end{rema}

\section{\sf Hochschild-Mitchell homology}\label{HM homology}

\subsection{\sf $G$-$k$-categories and graded $k$-categories }\label{section Gk and graded}

We begin this section by recalling the definition of the Hochschild-Mitchell homology of a $k$-category, see for instance \cite{mitchell, keller,mccarthy, cr}. Next we give properties of it for $G$-$k$-categories.

Moreover we provide the decomposition of the homology of a graded category along the conjugacy classes of $G$.  This will be used in Subsection \ref{section homology skew} for skew categories.

\Bdf
Let $\C$ be a $k$-category and let $C_{\bullet}(\C)$ be the chain complex given by:

$$C_n(\C) =\bigoplus_{x_0,x_1, \dots  ,x_n \in \C_{0}} \ _{x_0}\C_{x_n}\otimes\cdots\otimes\  _{x_2}\C_{x_1}\otimes\ _{x_1}\C_{x_0},$$
with  boundary map $d$ given by the usual formulas used to compute the Hochschild homology of an algebra, see for instance
\cite{loday, weibel, ce}.  Note that $C_0(\C) = \bigoplus_{x\in \C_0} \  _x\C_x$.

The \emph{Hochschild-Mitchell homology} $HH_*(\C)$ of $\C$ is the homology of the complex above.
\Edf

\subsubsection{\sf $G$-$k$-categories}

\begin{prop}
Let $\C$ be a $G$-$k$-category. The group $G$ acts on the chain complex $C_{\bullet}(\C)$ by automorphims.
\end{prop}
\begin{proo}
 For $s\in G$ and $(f_n\otimes \dots\otimes f_1\otimes f_0) \in C_n(\C)$ we define $$sf = (sf_n\otimes\dots\otimes sf_1\otimes sf_0)$$ and we have $d(sf)=sd(f)$. For simplicity we provide the verification for $n=2$, the general case follows the same kind of computations.
$$\begin{array}{llll}
d(sf)= &_{sx_0}(sf_2\ sf_1){}_{sx_1}\otimes {}_{sx_1}{(sf_0)}_{sx_0} -\   _{sx_0}\! (sf_2){}_{sx_2}\otimes \ _{sx_2}{(sf_1\ sf_0)}_{sx_0} +
\\&{}_{sx_1}{(sf_0\ sf_2)}_{sx_2}\otimes{}_{sx_2}{(sf_1)}_{sx_1}
\end{array}
$$
and
$$\begin{array}{llll}
sd(f) = &{}_{sx_0}\left( s(f_2f_1)\right)_{sx_1}\ \otimes {}_{sx_1}{(sf_0)}_{sx_0}  - \  _{sx_0}\! (sf_2){}_{sx_2} \otimes \ _{sx_2} \left(s(f_1f_0)\right)\! _{sx_0} \ +
\\& _{sx_1}\left(s(f_0f_2) \right)_{sx_2}\ \otimes \ _{sx_2} f_1{}_{x_1}.
\end{array}$$
\qed
\end{proo}
\Bcr
Let $\C$ be a $G$-$k$-category. The Hochschild-Mitchell homology $HH_* (\C)$ is a $kG$-module.
\Ecr
Recall that for a $kG$-module $M$, the $kG$-module of coinvariants of $M$ is
$$M_G = M/<sm-m>$$
where the denominator is the sub $kG$-module of $M$ generated by $\{sm -m \mid m\in M, s\in G\}.$ The module of coinvariants is the largest quotient of $M$ with trivial action of $G$. Considering $M$ as a $kG$-bimodule  with trivial action on the right, we have  $M_G= kG \otimes_{kG\otimes (kG)^{op}} M.$

If $G$ is of finite order invertible in $k$, then the coinvariants functor is exact. Moreover, $M_G$  is canonically isomorphic to the invariants
$$M^G= \{m\in M \mid sm =m \text{ for all } s \in G\}$$ through the morphism $m\mapsto \frac{1}{\vert G \vert} \sum_{s\in G}sm$.

\subsubsection{\sf Graded categories}

\Bdf
Let $G$ be a group.  A $k$-category $\B$ is \emph{$G$-graded} if for all $x,y \in \B_0$ there is a direct sum decomposition of vector spaces
 $$_y\B_x =\bigoplus_{s\in G} {}_y\B^s_x$$
 such that $_z\B_y^t \ _y\B^s_x\subset {}_z\B_x^{ts}$ for all objects $x,y,z\in \B_0$ and for every $s,t\in G.$ A morphism $f\in{}_y\B^s_x$ is called \emph{homogeneous of degree $s$} from $x$ to $y$, we often write it ${}_yf_x^s$ instead of $f$.
 \Edf

Next we provide a decomposition of the chain complex of a graded category. This corresponds to M. Lorenz \cite{lorenz} decomposition for the Hochschild homology of a $G$-graded  $k$-algebra (see also \cite{stefan1996}).

\begin{prop} Let $\B$ be a $G$-graded $k$-category. Let $D$ be a conjugacy class of $G.$ The following is a subcomplex of $C_{\bullet}(\B)$:
$$C^D_n(\B)=\bigoplus_{\substack{s_n \dots s_1s_0 \in D\\ x_n,\dots ,x_1,x_0 \in \B_0}} {}  _{x_0}\B^{s_n}_{x_n}\otimes \cdots
\otimes{}_{x_2}\B^{s_1}_{x_1}\otimes{}_{x_1}\B^{s_0}_{x_0}.$$
Let $Cl_G$ be the set of conjugacy classes of $G.$ There is a decomposition
$$C_{\bullet}(\B) = \bigoplus_{D\in Cl_G} C^D_\bullet (\B).$$

\end{prop}

\begin{proo}
We provide the verification for $n=2$, since it already gives the whole track to prove the result for all $n$. This check has the advantage of an easier reading, it presents the main ideas of the calculations, without too many technical details.

Let ${}_{x_0}f_2^{s_2}{}_{x_2} \otimes {}_{x_2}f_1^{s_1}{}_{x_1}\otimes{}_{x_1}f_0^{s_0}{}_{x_0} \in C_2^D(\B),$ with $s_2 s_1 s_0 \in D.$
We have
$$\begin{array}{cccccc}
d(f_2 \otimes f_1\otimes f_0)=\\f_2 f_1\otimes f_0  &-&f_2\otimes f_1 f_0 &+& f_0f_2\otimes f_1 \in\\

\left({}_{x_0} \B^{s_2s_1}_{x_1}\otimes {}_{x_1}\B^{s_0}_{x_0}\right)
&\oplus &\left(_{x_0}\B^{s_2}_{x_2}\otimes {}_{x_2}\B^{s_1s_0}_{x_0}\right)&\oplus
&\left({}_{x_1}\B^{s_0s_2}_{x_2}\otimes {}_{x_2}\B^{s_1}_{x_1}\right).
\end{array}
$$

Note that for the last summand $s_0s_2s_1 = s_0 (s_2 s_1 s_0)s^{-1}_0 \in D.$
\qed
\end{proo}
\Bte\label{decomposition HH of a graded}
Let $\B$ be a $G$-graded $k$-category and let $HH^D_*(\B)=H_*(C^D_\bullet(\B))$. There is a decomposition
$$ HH_*(\B) = \bigoplus_{D \in Cl_G} HH^D_*(\B).$$
\Ete

\subsection{\sf Skew categories}\label{section homology skew}

In this section we compare the coinvariants of the Hochschild-Mitchell homology of a $G$-$k$-category $\C$ with the trivial component of the homology of the skew category $\C[G]$. Indeed $\C[G]$ (defined below) is graded, hence we will use Theorem \ref{decomposition HH of a graded}.

\begin{defi}\cite{cm} Let $\C$ be a $G$-$k$-category. The \emph{skew category} $\C[G]$ has the same set of objects than $\C$. Let ${_y\C[G]_x}^s={}_y\C_{sx}$. The morphisms of $\C[G]$ from $x$ to $y$ are
\begin{equation}\label{eq:1}
_y\C[G]_x = \bigoplus_{s \in G}{{}_y\C[G]_x}^s.
\end{equation}
The composition  is defined through adjusting the first morphism in order to make it possible to compose it in $\C$ with the second one, as follows. If
$$_yf_{sx}\in{}_y\C_{sx} \subseteq \  _y\C[G]_x \mbox{ and } _zg_{ty} \in{}_z\C_{ty}\subseteq \  _z\C[G]_y, \mbox{then}$$
  $$(_zg_{ty})(_yf_{sx}) = {}_z(g\ \circ \ t\!f)_{tsx}\in {}_z\C[G]_x,$$
  where  $\circ$ denotes the composition of $\C$.
 \end{defi}

 \begin{rema}\label{definition skew group algebra}
 \

 \begin{enumerate}
  \item By definition, the direct summands of (\ref{eq:1}) are in one to one correspondence with elements of $G$. However some of the summands may be zero.
 \item  If $\C$ is a single object $G$-$k$-category with endomorphism algebra $\Lambda$, it is shown in \cite{cm} that the endomorphism algebra of the single object $k$-category  $\C[G]$ is the usual skew group algebra.  Namely as a vector space, $\Lambda[G]= \Lambda\otimes kG$. For $a,b\in \Lambda$ and $s,t\in G$, the product is given by
$$(a\otimes s)(b\otimes t) = as(b)\otimes st.$$

 \end{enumerate}
 \end{rema}

 \begin{lemm}
 The category $\C[G]$ is graded.
 \end{lemm}
 \begin{proo}
 We have
  \begin{equation}\label{observe}
  \ {{}_z\C[G]_{y}}^t  \  {{}_y\C[G]_x}^s \subset\  {{}_z\C[G]_x}^{ts}.
 \end{equation}
 \qed
    \end{proo}

Hence the decomposition of Theorem \ref{decomposition HH of a graded} holds for the homology of $\C[G]$.

\subsubsection{\sf Free action}

To prove the next result  we first need to restrict ourselves to a free action. Later we will be able to consider the general case by using the resolving category of Section \ref{section resolving}.

\Bte\label{homology free}
Let $\C$ be a $G$-$k$-category with free action of $G$. Let $\{1\}$ be the trivial conjugacy class of $G$. There is an isomorphism

$$HH^{\{1\}}_* (\C[G]) \simeq  H_* \left((C_\bullet (\C)_G\right).$$
If the coinvariants functor is exact then
$$HH^{\{1\}}_* (\C[G]) \simeq  H_* (\C)_G.$$
\Ete

Next we provide properties of a skew category that we will need for proving Theorem \ref{homology free}. The following result is proved in \cite{cm}, we give a proof for completeness.

\Blm \label{isos}
Let $\C$ be a $G$-$k$-category. Let $x$ and $y$ be objects in the same orbit of the action of $G$. They are isomorphic in $\C[G].$
\Elm

\Bdm
Let $t\in G$ such that $y=tx$. Recall that $_{tx}\C[G]_x = \bigoplus_{s\in G} {}_{tx}\C_{sx}.$ Let
$$a = {}_{tx}1_{tx}\in{}_{tx}\C[G]^t_x= {}_{tx}\C_{tx}  \mbox{ \ and\ }  b= {}_x1_x \in{}_x\C[G]^{t^{-1}}_{tx} = {}_{x}\C_{x}.$$
Using the composition defined in $\C[G]$, we obtain that $a$ and $b$ are mutual inverses in the skew category.\qed
\Edm
\Bdf
Let $G$ be a group acting on a set $E.$ A \emph{transversal $T$ of the action} is a subset of $E$ consisting of exactly one element in each
orbit of the action.
\Edf
Equivalently $T\subset E$ is a transversal if for each $x\in E$ there exists a unique $u(x) \in T$ such that there exists some $s\in G$ with $ x =su(x)$.

Note that the action is free if and only if for each $x\in E$, there exists a unique $s\in G$ such that $x=su(x)$.
\Blm \label{lemma transversal}
Let $\C$ be a $G$-$k$-category, let $T\subseteq \C_0$ be a transversal of the action of $G$ on $\C_0$, and let $\C_T[G]$ be
the full subcategory of $\C[G]$ with set of objects $T.$ For each conjugacy class $D$ of $G$ we have
$$HH^D_*(\C_T[G]) =  H_*^D (\C[G]).$$
\Elm

\Bdm  In \cite{herscovich} E. Herscovich (see also \cite{bw,hsja}) proved that if $\C$ and $\D$ are $k$-categories and if $F:\C\to \D$ is a $k$-equivalence, then $F$ induces a quasi-isomorphism $C_\bullet(\C)\to C_\bullet(\D)$. Moreover, if $\C$ and $\D$ are $G$-graded and if $F$ is homogeneous, then the induced quasi-isomorphism clearly preserves the decomposition along conjugacy classes of $G$. Since $T$ is a transversal, the above Lemma \ref{isos} shows that the inclusion functor $\C_T[G]\subset \C[G]$ is dense. Moreover it corresponds to a full subcategory, hence it is full and faithful, in addition of being homogeneous. The induced quasi-isomorphism provides then the result. \qed
\Edm

\noindent\textbf{Proof of Theorem \ref{homology free}.}  As previously, we provide the proof for $n=2$. This avoids lengthy and uninformative computations.

Let $T$ be a transversal of the free action of $G$ on $\C_0$.
In order to define an isomorphism of chain complexes $A:C_\bullet (\C)_G \to C_\bullet (\C^{\{ 1\}}_{T}[G]),$  let
$_{x_0}f_2{}_{x_2} \otimes {}_{x_2}f_1{}_{x_1} \otimes {}_{x_1} f_0{}_{x_0}$ be a chain  of $C_2 (\C).$
In $\left(C_2(\C)\right)_G,$ we begin by modifying the chain in order that the starting (and hence the ending)
object $x_0$ belongs to $T.$ More precisely, there exists a unique $s\in G$ such that $sx_0 ={} u_0 \in T.$ Then
$$f_2\otimes f_1 \otimes f_0 \equiv s(f_2\otimes f_1 \otimes f_0)=\
_{u_0}\!(sf_2){}_{sx_2} \otimes {}_{sx_2}(sf_1){}_{sx_1} \otimes {}_{sx_1} (sf_0){}_{u_0}.$$
In other words we can assume that the chain is of the form
$$_{u{}_0}(f_2){}_{x_2} \otimes {}_{x_2}(f_1){}_{x_1} \otimes {}_{x_1} (f_0){}_{u_0}\text { for } u_0 \in T.$$
For $i=1,2$, let $u_i = u(x_i)$ be the unique element of $T$ which is in the orbit of $x_i$.  Moreover let $s_i$ be the unique element of $G$ such that $x_i = s_iu_i$. We define
$$\begin{array}{llll}
A\left(_{u{}_0}(f_2){}_{s_2u_2}\ \otimes \ {}_{s_2 u_2}(f_1){}_{s_1u_1}\   \otimes\ {}_{s_1u_1} (f_0){}_{u_0}\right) =& \\
{}_{u_0}( f_2) {}_{s_2 u_2}\ \otimes   \  _{u_2}(  s_2^{-1}f_1)_{s_2^{-1}s_1u_1}\ \otimes\  _{u_{1}}(s^{-1}_1f_0)_{s_1^{-1} u_0} =&\\
_{u_0}[f_2]^{s_2}_{u_2}\ \otimes{}_{u_2}[s^{-1}_2 f_1]_{u_1}^{s_2^{-1}s_1}\ \otimes{}_{u_1}[s_1^{-1} f_0]^{s_1^{-1}}_{u_0}.&
\end{array}
$$
This chain belongs to $C_2^{\{1\}}(\C_T[G])$ since $s_2(s_2^{-1}s_1)s_1^{-1} =1.$ For a 3-chain the formula defining $A$ is
$$\begin{array}{llll}
A(_{u{}_0}f_3{}_{s_3u_3} \otimes {}_{s_3u_3}
{f_2} _{s_2u_2} \otimes  {}_{s_2 u_2}f_1{}_{s_1u_1}  \otimes {}_{s_1u_1} f_0 {}_{u_0}) = & \\
{}_{u_0}[f_3]^{s_3}_{u_3}\ \otimes \
_{u_3} [s^{-1}_3( f_2)]  _{u_2}^{ s_3^{-1}s_2}\ \otimes \ {}_{u_{2}}{[(s^{-1}_2f_1)]}^{s_2^{-1} s_1}_{u_1}\ \otimes\ {}
 _{u_1}{[s_1^{-1} f_0]}^{{s_1}^{-1}}_{u_0}.  &
 \end{array}
$$
Next we verify that $A$ is a chain map.
$$
\begin{array}{lll}
dA(f_2\otimes f_1 \otimes f_0) &=&
{[f_2]^{s_2}[s^{-1}_2 f_1]^{s_2^{-1}s_1}}\otimes [s^{-1}_1 f_0]^{s_1^{-1}} -\\&&
[f_2]^{s_2}\otimes [s^{-1}_2 f_1]^{s_2^{-1}s_1}[s^{-1}_1 f_0]^{s_1^{-1}} + \\&&
[s^{-1}_1 f_0]^{s_1^{-1}} [f_2]^{s_2}\otimes [s^{-1}_2 f_1]^{s_2^{-1}s_1} \\\\&=&

[f_2f_1]^{s_1}\otimes [s_1^{-1} f_0]^{s_1^{-1}} - \\ &&[f_2]^{s_2} \otimes [(s_2^{-1}f_1) (s^{-1}_2 f_0)]^{s_2^{-1}} + \\&& {[(s_1^{-1}f_0) (s_1 ^{-1} f_2)]^{s_1^{-1}s_2} \otimes
[s_2^{-1} f_1]^{s_2^{-1}s_1}.}
\end{array}
$$
Recall that
$$
\begin{array}{lll}
d(f_2\otimes f_1 \otimes f_0) =  &{}_{u_0}(f_2f_1)_{s_1u_1}\otimes {}_{s_1u_1}(f_0)_{u_0} - \\& {}_{u_0}(f_2)_{s_2u_2}\otimes {}_{s_2u_2}(f_1f_0)_{u_0}
+\\& {}_{s_1u_1}(f_0f_2)_{s_2u_2}\otimes{}_{s_2u_2}(f_1)_{s_1u_1}.
\end{array}
$$
In order to compute $Ad$, notice that up to the action, that is in the coinvariants, the last term of the previous sum can be rewritten:
$$ {}_{s_1u_1}(f_0f_2)_{s_2u_2}\otimes{}_{s_2u_2}(f_1)_{s_1u_1}\equiv \  {}_{u_1}(s^{-1}(f_0f_2))_{s_1^{-1}s_2u_2}\otimes{}_{s_1^{-1}s_2u_2}(s^{-1}_1f_1)_{u_1}.$$

 This way the last summand of $d(f_2\otimes f_1 \otimes f_0)$ starts and ends at $u_1 \in T$, which is required in order to apply $A.$ Hence
 $$
\begin{array}{llllllllll}
 Ad(f_2\otimes f_1 \otimes f_0) &=&A\big({}_{u_0}(f_2f_1)_{s_1u_1}\otimes {}_{s_1u_1}(f_0)_{u_0} -
  {}_{u_0}(f_2)_{s_2u_2}\otimes{}_{s_2u_2}(f_1f_0)_{u_0} + \\
  &&{}_{u_1}(s_1^{-1}(f_0f_2))_{s^{-1}_1s_2u_2}\otimes \ _{s^{-1}_1s_2u_2}(s_1^{-1}f_1)_{u_1}\big)\\\\&=&
[(f_2f_1)]^{s_1} \otimes [s_1^{-1}f_0]^{s_1^{-1}}  - [f_2]^{s_2}\otimes [s_2^{-1}(f_1f_0)]^{s_2^{-1}} +\\
&&[s_1^{-1}(f_0f_2)]^{s_1^{-1}s_2}\otimes [(s_1^{-1}s_2)^{-1} s_1^{-1}f_1]^{(s_1^{-1}s_2)^{-1}}
\end{array}
$$
shows that $Ad =dA.$

Let $g_2\otimes g_1\otimes g_0 \in C_2^{\{1\}}(\C_T[G]),$  that is
 $$g_2\otimes g_1\otimes g_0= {}^{}_{u_0}[g_2]^{s_2}_{u_2} \otimes{}_{u_2}[g_1]^{s_1}_{u_1} \otimes{}_{u_1}[g_0]^{s_0}_{u_0} = {}_{u_0}g_{2}{}_{s_2u_2} \otimes{}_{u_2}g_1{}_{s_1u_1} \otimes {}_{u_1}g_0{}_{s_0 u_0}$$
 with $s_2s_1s_0 =1.$  Let

 $$B : C_{\bullet}^{\{1\}} (\C_T[G]) \to \left(C_\bullet(\C)\right)_G$$
 be defined by
 $$B (g_2\otimes g_1\otimes g_0) = {}_{u_0}(g_2){}_{s_2u_2}\otimes {}_{s_2u_2}(s_2g_1)_{s_2s_1u_1}\otimes {}_{s_2s_1u_1}(s_2s_1 g_0)_{s_2s_1s_0u_0}.$$
We observe that since $s_2s_1s_0 =1,$ we have that $s_2s_1s_0u_0 = u_0$.
Next we will show that $A$ and $B$ are mutual inverses. This will imply that $B$ is a chain map, since $A$ is a chain map.
$$\begin{array}{llll}
AB(g_2\otimes g_1\otimes g_0) = \\ A(_{u_0}g_2{}_{s_2u_2}\otimes {}_{s_2u_2}(s_2g_1)_{s_2s_1u_1}\otimes {}_{s_2s_1u_1}(s_2s_1 g_0)_{u_0}) = \\
{}_{u_0}g_2{}_{s_2u_2}\otimes {}_{u_2}(s_2^{-1}s_2g_1)_{s_1u_1}\otimes {}_{u_1}((s_2s_1)^{-1}(s_2s_1) g_0)_{(s_2s_1)^{ -1}u_0}.
\end{array}
$$

Since $s_2 s_1 s_0 =1$, we obtain
$$_{u_0}g_2{}_{s_2u_2}\otimes {}_{u_2}g_1{}_{s_1u_1}\otimes {}_{u_1}g_0{}_{s_0u_0} =
{}^{}_{u_0}[g_2]^{s_2}_{u_2} \otimes{}^{}_{u_2}[g_1]^{s_1}_{u_1} \otimes{}^{}_{u_1}[g_0]^{s_0}_{u_0}. $$
$$\begin{array}{lll}
 BA\left(_{u{}_0}f_2{}_{s_2u_2} \otimes {}_{s_2 u_2}f_1{}_{s_1u_1}  \otimes {}_{s_1u_1} f_0{}_{u_0}\right)=\\
 B\left(_{u{}_0}[f_2]{}_{u_2}^{s_2} \otimes {}_{ u_2}[s_2^{-1}f_1]^{s_2^{-1}s_1} _{u_1} \otimes {}_{u_1} [s_1 ^{-1}f_0)]^{s_1^{-1}} _{u_0}\right)
 =\\ {}_{u_0}f_2{}_{s_2u_2} \otimes{}_{s_2u_2}(s_2s_2^{-1}f_1)_{s_1u_1} \otimes{}_{s_2s_2^{-1}s_1u_1} ((s_2s_2^{-1}s_1)s_1^{-1}f_0)_{(s_2s_2^{-1}s_1)s_1^{-1} u_0}
= \\ {}_{u{}_0}f_2{}_{s_2u_2} \otimes {}_{s_2 u_2}f_1{}_{s_1u_1}  \otimes {}_{s_1u_1} f_0{}_{u_0}.
\end{array}$$
\qed

\subsubsection{\sf General case}\label{subsection general case homology}

 The next Lemma enables to generalise Theorem \ref{homology free} to a $G$-$k$-category where the action of $G$ is non necessarily  free. Recall that $M_G(\C)$ is the resolving category of a $G$-$k$-category, see Definition \ref{MGC}.

\Blm Let $\C$ be a $G$-$k$-category. The chain complexes $C_\bullet (M_G(\C))$ and $C_\bullet (\C)$ are $kG$-quasi-isomorphic.
\Elm

\begin{proo}
By Theorem \ref{L}, there is an equivalence of categories $L:M_G(\C) \to \C$. As for Lemma \ref{lemma transversal}, we infer from \cite{herscovich}  that there is an induced quasi-isomorphism
$\C_\bullet (M_G(\C)) \to C_\bullet (\C).$
\qed
\end{proo}

\Bte
Let $\C$ be a $G$-$k$-category.
$$HH_*^{\{1\}}(\C[G]) = H_*\left(C_\bullet (\C)_G\right).$$
If the coinvariants functor is exact
$$HH_*^{\{1\}}(\C[G]) = HH_*(\C)_G.$$
\Ete

\Bdm Due to Theorem \ref{homology free} the result holds for the $G$-$k$-category $M_G(\C).$
The equivalence of $G$-$k$-categories $L:M_G(\C) \to \C$ provides a homogeneous equivalence of $G$-graded $k$-categories
$$L[G]: M_G(\C)[G] \to \C[G]$$ which gives a quasi-isomorphism
$$C_\bullet (M_G(\C)[G]) \to C_\bullet (\C[G])$$
preserving the decomposition of chain complexes along the conjugacy classes of $G$. Hence
$$HH^{\{1\}}_*(M_G(\C)[G]) = HH^{\{1\}}_*  (\C[G]).$$
By the above proposition
$$H_*( C_\bullet(M_G(\C))_G) = H_* (C_\bullet(\C)_G).$$
If the coinvariants functor is exact then
$$HH_*\left(M_G(\C)\right)_G = \left(HH_* (\C)\right)_G.$$
\qed
\Edm

\subsection{\sf Galois coverings}\label{section Galois}

In this section we reformulate Theorem \ref{homology free} in terms of Galois coverings. First we recall the definition of a
quotient category.
\Bdf (see \cite{bg,ri} and also \cite{cm})\label{definition Galois}
Let $\C$ be a $G$-$k$-category with a free action of $G$. The \emph{quotient category} $\C/G$ has set of
objects the set of orbits $\C_0/G.$ Let $\alpha$ and  $\beta$ be orbits. The vector space of morphisms from $\alpha$ to $\beta$  is
$$ _\beta(\C/G)_\alpha = \left(\bigoplus_{\substack{x\in \alpha \\ y \in \beta}}{}_y\C_x\right)_G.$$
 Let $\gamma, \beta, \alpha$ be orbits. Let
$$g\in{}_z\C_{y'}\text{ and } f\in{}_y\C_x,$$ where $z\in \gamma$, $y \text{ and } y' \in \beta$, and  $x\in\alpha.$
Let $s$ be the unique element of $G$ such that $sy = y',$ then $f\equiv sf$ in the coinvariants. The composition $gf$ in $\C/G$ is
 $$g f ={}_{z}g_{y'} \ {}_{sy}sf_{sx}\in{}_{\gamma}(\C/G)_\alpha.$$
 There is no difficulty in verifying that this is a well defined associative composition.

\Edf

\Bdf
A Galois covering of $k$-categories is a functor $\C\to \C/G,$ where $\C$ is a $G$-$k$-category with free action,  and where the functor is the canonical projection functor.
\Edf

If $\C_0$ is finite we infer from a Galois covering as above a map $a(\C)\to a(\C/G)$. Note that $a(\C)$ is not an algebra map in general, see Remark \ref{notalgebramaps}.

Let $\C\to \C/G$ be a Galois covering and let $T$ be a transversal of the action of $G$ on $\C_0$. For each
orbit $\alpha$, let $u_{\alpha}\in T$ be the unique element of $T$ which belongs to $\alpha.$

It is shown in Lemma 2.2 of \cite{cm} that through a canonical identification we have
$$_\beta (\C/G)_\alpha = \bigoplus_{ s\in G}{}_{u_\beta}\C_{su_\alpha}.$$
This provides that $\C/G$ is graded by $G$. Indeed let ${}_{\beta}^{}(\C/G)^s_{\alpha}={}_{u_\beta}\C_{su_\alpha}$, and notice that
$$(_{u_\gamma}g_{tu_\beta})({}_{u{}_{\beta}}f_{su_\alpha}) = \left(_{u_\gamma}g_{tu_\beta})({}_{tu{}_{\beta}}tf_{tsu_\alpha}\right)=
 {}_{u_\gamma}(g  (tf))_{tsu_\alpha} \in {}_\gamma^{} (\C/G)_{\alpha}^{ts}{} .$$
The following result can be deduced from \cite{cm}. We provide a proof for completeness.
\Bpo\label{Galois and C_T are homogeneous equivalent}
Let $\C\to \C/G$ be a Galois covering. Let $T\subset \C_0$ be a transversal of the free action of $G$, and consider the
$G$-grading of $\C/G$ determined by $T.$

The graded categories $\C/G$ and $\C_T[G]$ are isomorphic by a homogeneous functor, that is by a functor which is homogeneous on morphisms and graded with respect to composition of morphisms.
\Epo

\Bdm
There is a bijection between the objects of $\C/G$ and those of $\C_T[G]$. The previous considerations shows that the morphisms of both
categories are identified in a homogeneous manner. Moreover we have already used that the inclusion $C_T[G] \subset \C[G]$ provides an equivalence of categories, and this equivalence is homogeneous.
\qed
\Edm

The above analysis provides the translation of  Theorem \ref{homology free} in terms of Galois coverings, as follows.
\Bte \label{theorem cohomology Galois covering}
Let $\C \to \C/G$ be a Galois covering.
$$HH_*^{\{1\}}(\C/G) = H_*(C_\bullet(\C)_G).$$
If the coinvariants functors is exact
$$HH_*^{\{1\}}(\C/G) = HH_*(\C)_G.$$
\Ete

\section{\sf Hochschild-Mitchell cohomology}\label{HM cohomology}

\subsection{\sf Graded categories}\label{section cohomology graded}

As mentioned in the Introduction, Hochschild-Mitchell cohomology is more intricate than homology since it makes use of direct products. We begin this section by recalling its definition. If the category is graded, we provide the decomposition along conjugacy classes.

\Bdf Let $\C$ be a $k$-category. Let $C^\bullet (\C)$ be the complex of cochains given by:
$$\C^0(\C) = \prod_{x\in \C_0}{}_x\C_x,$$
$$ \C^n(\C) = \prod_{x_{n+1}, \dots , x_1 \in \C_0} C_{x_{n+1}, \dots, x_1}  \text{ for } n> 0$$
where
$$C_{x_{n+1}, \dots, x_1} = \mathsf{Hom}_k({}_{x_{n+1}}\C_{x_n}\otimes \cdots \otimes {}_{x_2}\C_{x_1},\ {}_{x_{n+1}}\C_{x_1} )
$$
The coboundary $d$ is given by the formulas below which are the usual ones for computing Hochschild cohomology, see for instance \cite{ce,weibel}.

Let $\varphi$ be a cochain of degree $n$, that is  a family  of $k$-morphisms $\varphi =\{ \varphi_{(x_{n+1}, \dots, x_1)}\}$. Its coboundary $d\varphi$ is the family
$\{ (d\varphi)_{(x_{n+2}, \dots , x_1)} \} $
given by
\begin{equation}\label{d}
\begin{array}{lllll}
(d\varphi)_{(x_{n+2}, \dots, x_1)}\left({}_{x_{n+2}}(f_{n+1})_{x_{n+1}}\otimes \cdots \otimes {}_{x_2}(f_1)_{x_1}\right)=\\
(-1)^{n+1} f_{n+1} \varphi_{(x_{n+1},\dots,x_1)}(f_n\otimes \cdots\otimes f_1)
 +\\
 \sum_{i=1}^n (-1)^{i+1} \varphi_{(x_{n+1},\dots, x_{i+2}, x_i, \dots x_1)}(f_{n+1}\otimes\cdots\otimes f_{i+1}f_i\otimes \cdots\otimes f_1) + \\
\varphi_{(x_{n+2}, \dots x_2)}(f_{n+1}\otimes \cdots  \otimes f_2) f_1.
\end{array}
\end{equation}

Note that $d$ is well defined since the spaces of cochains are direct
products. The\emph{ Hochschild-Mitchell cohomology }of $\C$ is $HH^*(\C) = H^*(C^\bullet(\C)).$
\Edf
\begin{rema}
The compatibility with the Hochschild cohomology of a $k$-algebra arises from the following.
Let $\C$ and $\D$ be $k$-categories,  let $\C\otimes_k \B$ be their tensor product and let $\mathsf{Fun}_k(\C,\D)$ be the $k$-functors from $\C$ to $\D$.
Then for $n\geq 0$ we have
$$\C^n(\C)= \mathsf{Fun}_k (\C^{\otimes_kn}, \C).$$
Note that $\C^{\otimes_k0}$ is the following category. The objects are the same than those of $\C$. The space of morphisms between any pair of different objects is zero. The endomorphisms of each object is the $k$-algebra $k$.
\end{rema}

In  degree zero we have
$$ HH^0(\C) = \{{}(_xf_x)_{x\in \C_0} \mid {}_yg_x{}\ _xf_x = {}_yf_y{}\ _yg_x   \text{ for all }{}  _yg_x \in {}_y\C_x\}.
$$
As for Hochschild cohomology of algebras, the cup product is defined  at the cochain level as follows. Let $\varphi \in  \C_{x_{n+1}, \dots, x_1}$ and $ \psi \in \C_{y_{m+1}, \dots, y_1}$. If $x_{n+1} \neq y_1$ the cup product $ \psi \smile \varphi$ is zero.
Otherwise the cup product  $\psi \smile \varphi \in \C_{y_{m+1}, \dots, y_1, x_n,\dots, x_1}$ is
$$(\psi \smile \varphi)(f_{n+m}\otimes\cdots\otimes f_1) = \psi(f_{n+m}\otimes\cdots \otimes f_{n+1}) \varphi(f_n\otimes\cdots\otimes f_1).$$
The cup product verifies the graded Leibniz rule, and it provides a graded commutative $k$-algebra structure on $HH^*(\C).$ In particular
$HH^0(\C)$ is a commutative $k$-algebra which is the center of the category.

\Bpo
Let $\B$ be a $G$-graded category, and let $Cl(G)$ be the set of conjugacy classes of $G.$ There is a decomposition
$$HH^*(\B) = \prod_{D\in Cl(G)}HH^*_D(\B)$$ where $HH_{\{1\}}^*(\B)$ is a subalgebra.
\Epo
\Bdm
For $D\in Cl(G)$ we  provide a subcomplex of cochains $C^\bullet_D(\B) $ of $\C^\bullet(\B)$ as follows.
Let $\varphi$ be a cochain of degree $n$. We say that $\varphi$ is\emph{ homogeneous of type $(s_n,\dots,s_1, s_0)$}  if:
\begin{enumerate}
 \item Each component of $\varphi$ has its image contained in the set of homogeneous morphisms of degree $s_0$.
\item For  $(s'_n,\dots,s'_1)\neq (s_n,\dots,s_1)$, each component of $\varphi$ restricted to tensors of homogeneous morphisms degree $(s'_n,\dots,s'_1)$
 is zero.
\end{enumerate}
The formula (\ref{d}) which defines the coboundary $d$ has $n+2$ summands. Let $d_{n+1}$ be the first one, let $d_0$ be the
last one, and let $d_i$ denote the in between summands indexed according to the appearance of the composition
$ ``f_{i+1} f_i" \text{ for } i=n,\dots, 1.$

Let $\varphi$ be homogeneous of type $(s_n,\cdots,s_1, s_0).$ We observe the following:

\begin{itemize}
\item $d_{n+1}\varphi$ is a sum of homogeneous cochains of types $(s, s_n, \dots, s_1, ss_0)$  for $s\in G.$
\item $d_i\varphi$ is a sum of homogeneous cochains of types
$$ (s_n, \dots , s_{i+1}, s'', s', s_{i-1}, \dots, s_1, s_0)$$ for
$s'', s'\in G$ with $s''s'= s_i. $
\item $d_0\varphi$ is a sum of homogeneous cochains of types $(s_n, \dots, s_1, s, s_0 s)$ for $s\in G.$
\end{itemize}
Let us call the \emph{class} of $(s_n, \dots, s_1, s_0)$  the product $s_n\dots s_1s_0^{-1}\in G.$
The above considerations show that if $\varphi$ is homogeneous of type $(s_n, \dots, s_1,s_0)$, that is of class
$c= s_n\dots s_1 s_0^{-1}$, then $d\varphi$ is a sum of homogeneous cochains. Note that these cochains are in general of different sorts, but their classes are all conjugated to $c.$

Let $D$ be a conjugacy class and let $C^\bullet_D(\C)$ be the homogeneous cochains which classes of types are in $D.$ We have showed that $C^\bullet_D(\C)$ is a cochains subcomplex of $C^\bullet (\C).$  Moreover
$$ C^\bullet(\C) = \prod_{D \in Cl(G)} C^\bullet_D(\C).$$
Clearly, if $\varphi \text{ and } \psi$ are homogeneous cochains which classes of types are both 1, then $ \psi \smile \varphi $ is also of class type 1.
\Edm
\qed

\subsection{\sf Skew categories}

We recall that the skew category is graded, so results of the Subsection \ref{section cohomology graded} are in force. As for homology, at a first glance we are only able to provide a result if the action is free. The resolving category of Section \ref{section resolving} enables to extend the result for any $G$-$k$-category.

Firstly we provide some  properties of the cohomology of a $G$-$k$-category that we need in the sequel.

If $\C$ is a $G$-$k$-category then $C^n (\C)$ is a $kG$-module as follows. For $n>0$, let
$\varphi =\{ \varphi_{(x_{n+1}, \dots, x_1) }\}\in C^n(\C)$ be a cochain, where
 $$\varphi_{(x_{n+1}, \dots, x_1)}: {}_{x_{n+1}}\C_{x_n}\otimes\cdots\otimes{}_{x_2}\C_{x_1} \to{}_{x_{n+1}}\C_{x_1}.$$
 Let $s\in G$  and let
 $$s.[\varphi_{(x_{n+1}, \dots, x_1)}]: {}_{sx_{n+1}}\C_{sx_n}\otimes\cdots\otimes{}_{sx_2}\C_{sx_1} \to{}_{sx_{n+1}}\C_{sx_1}$$
be defined by
$$s.[\varphi_{(x_{n+1},\dots,x_1)}](f_n\otimes\cdots\otimes f_1) = s[\varphi_{(x_{n+1},\dots,x_
1)}(s^{-1}f_n\otimes\cdots\otimes s^{-1}f_1)] .$$
We set  $s.\varphi = \{ s.[\varphi_{x_{n+1}, \dots, x_1}]\}.$

For $n=0$, let $\varphi =\{ \varphi_{x}\}\in C^0(\C)$, where $\varphi_{x}\in {}_x\C_x.$ For $s\in G$ we have $s\varphi_x\in{}_{sx}\C_{sx}$. We set $s.\varphi = \{ s\varphi_{x}\}$.

\Brm
Let $(\ )^G$ be the invariants functor. Then  $\varphi\in (C^n(\C))^G$ if and only if for all $s\in G$ and for any sequence
of objects $x_{n+1}, \dots, x_1$ we have that
$$\varphi_{(sx_{n+1},\dots,sx_1)}(sf_n\otimes\cdots\otimes sf_1) =s\left[\varphi_{(x_{n+1},\dots,x_1)}(f_n\otimes\cdots\otimes f_1)\right].$$
\Erm
Clearly the action of $G$ commutes with the coboundary of $C^\bullet(\C)$. Moreover the action of $G$ is by automorphisms of the algebra structure given by the cup product. In other words $C^\bullet (\C)$ is a graded differential algebra with an action of $G$ by automorphisms of its structure.

In particular $\left(C^\bullet(\C)\right)^G$ is a graded differential algebra. Moreover the inferred action of $G$ on $HH^*(\C)$ is by automorphisms of the algebra. If the invariants functor is exact, then $\left(HH^*(\C)\right)^G = H^*(C^\bullet (\C)^G)$ as
$k$-algebras.

\subsubsection{\sf Free action}
\Bte\label{cohomology free}
Let $\C$ be a $G$-$k$-category with a free action of $G$, and let $\C[G]$ be the $G$-graded skew category. There is an isomorphism of $k$-algebras

$$HH^*_{\{1\}}(\C[G])\simeq H^*\left(\left(C^\bullet(\C\right)^G\right).$$

 If the invariants functor is exact, we infer an isomorphism of $k$-algebras

 $$HH^*_{\{1\}}(\C[G])\simeq HH^*(C)^G.$$
\Ete

\Bdm
Let $T$ be a transversal of the action of $G$ on $\C_0$ and let $\C_T[G]$ be the full subcategory of $\C[G]$ with
set of objets $T.$ Let $D$ be a conjugacy class of $G.$ We assert that
$$HH^*_D(\C_T[G]) = HH^*_D(\C[G]).$$
Indeed the equivalence of categories given by the inclusion $\C_T[G] \subseteq \C[G]$ induces a quasi-isomorphism of the complexes of cochains which preserves the decomposition along the conjugacy classes of $G.$

Moreover for the trivial conjugacy class the quasi-isomorphism is a morphism of differential graded algebras. Hence
$$HH^*_{\{1\}}(\C_T[G]) = HH^*_{\{1\}}(\C[G])$$ as $k$-algebras.

In what follows we will prove that there are morphisms of graded differential algebras
$$\left(C^\bullet (\C)\right)^G\underset{B}{\stackrel{A}{\rightleftarrows} }C^\bullet_{\{1\}}(\C_T[G])$$ which are mutual inverses.

We just give the proof for $n=3$, which is easier to read than the general case. Moreover, the main idea of the calculation is already highlighted -- the general case would unnecessarily lengthen the text.

 Let $\psi \in \left(C^3(\C)\right)^G$ and let $u_4, u_3, u_2, u_1 \in T.$ We will define $(A\psi)_{(u_4, u_3, u_2, u_1)}$ on each homogeneous component.

 Let
  $$f_3\otimes f_2 \otimes f_1 \in{}^{}_{u_4}\C_T[G]^{s_3}_{u_3}\otimes{}^{}_{u_3}\C_T[G]^{s_2}_{u_2} \otimes{}^{}_{u_2}\C_T[G]^{s_1}_{u_1}. $$
Recall that by the definition of the morphisms of $\C[G]$ we have that $f_i \in{}_{u_{i+1}}\C_{t_iu_i}$ for $i=1,2,3.$
Let
$$(A\psi)_{(u_4, u_3, u_2, u_1 )}(f_3\otimes f_2 \otimes f_1) = \psi _{(u_4, s_3u_3, s_3s_2u_2, s_3s_2s_1u_1)}(f_3\ \otimes s_3 f_2\otimes s_3s_2 f_1).$$
We observe that this definition makes sense since
$$f_3\otimes s_3 f_2\otimes s_3s_2 f_1\in {}_{u_4}\C_{s_3u_3}\otimes{}_{s_3u_3}\C_{s_3s_2u_2}\otimes{}_{s_3s_2u_2}\C_{s_3s_2s_1u_1}.$$
Moreover $$A\psi (f_3\otimes f_2 \otimes f_1)\in {}_{u_4}\C_{s_3s_2s_1u_1} = {}^{}_{u_4}\C_T[G]^{s_3s_2s_1}_{u_1},$$
that is we have indeed defined a homogeneous cochain of type $(s_3, s_2, s_1, s_3s_2s_1)$, which is of class $\{1\}.$

The verification that $dA = Ad$ is straightforward, it uses in a crucial way that $\psi$ is an invariant; the  formulas defining the composition in $\C[G]$ are required as well . Analogously, it is easy to verify that $A(\psi '\smile \psi) = A(\psi ')\smile A(\psi)$.

Let $\varphi \in C^3_{\{1\}}(\C_T[G]).$ In order to define $(B\varphi)_{(x_4, x_3, x_2, x_1)}$ we first observe that since the action
of $G$ on $\C_0$ is free, there exist $s_4, s_3, s_2, s_1\in G $ which are unique such that $x_i=s_iu_i$ for $i=1,2,3,4$.

 Let
 $g_3\otimes g_2\otimes g_1\in{}_{s_4u_4}\C_{s_3u_3}\otimes{}_{s_3u_3}\C_{s_2u_2}\otimes{}_{s_2u_2}\C_{s_1u_1}.$
We define $(B\varphi)_{(x_4, x_3, x_2, x_1)}$ as follows:
 $$ (B\varphi)(g_3\otimes g_2\otimes g_1) = s_4\varphi_{(u_4, u_3, u_2, u_1)}(s_4^{-1}g_3 \otimes s_3^{-1}g_2\otimes s^{-1}_2g_1).$$
 In order to verify that this is well defined, note first that

 %\vskip3mm
 $$\begin{array}{llllllllllllllll}
s_4^{-1}g_3 \otimes s_3^{-1}g_2\otimes s^{-1}_2g_1 \in & {}_{u_4}\C_{s^{-1}_4 s_3 u_3} & \otimes &{}_{u_3}\C_{s^{-1}_3s_2u_2} &\otimes  &{}_{u_2}\C_{s_2^{-1}s_1u_1} =\\
 &{}_{u_4}\C[G]^{s^{-1}_4s_3}_{u_3} & \otimes &{}_{u_3}\C[G] ^{s_3^{-1}s_2}_{u_2} &\otimes &{}_{u_2}\C[G] ^{s_2^{-1}s_1}_{u_1}.
  \end{array}$$

 Secondly, using that $\varphi$ is a cochain for the trivial conjugacy class, we obtain
$$\begin{array}{lllllll}
\varphi_{(u_4, u_3, u_2, u_1)} (s^{-1}_4 g_3, s^{-1}_3 g_2, s_2^{-1}g_1) \in &{}_{u_4}\C[G] ^{s_4^{-1}s_3 s_3^{-1}s_2s_2^{-1}s_1}_{u_1} =& {}_{u_4}\C[G] ^{s_4^{-1}s_1}_{u_1}= \\ &&{}_{{}_{u_4}}\C_{s_4^{-1}s_1u_1}.
 \end{array}$$
 Hence $(B\varphi)(g_3\otimes g_2 \otimes g_1) \in {}_{su_4}\C_{su_1}$, therefore $B\varphi \in C^3(\C).$ Next we check that $B\varphi$ is an invariant cochain. Let $t\in G$, we assert that
 $$\begin{array}{ll}
  t(B\varphi)_{(s_4u_4, s_3u_3, s_2u_2, s_1u_1)}(g_3\otimes g_2\otimes g_1) = \\  B\varphi _{(ts_4u_4, ts_3u_3, ts_2u_2, ts_1u_1)}(tg_3\otimes tg_2\otimes tg_1).
 \end{array}$$
Indeed, the second term is by definition $$ ts_4\ \varphi_{(u_4, u_3, u_2, u_1)}\left((ts_4)^{-1} tg_3\otimes (ts_3)^{-1} tg_2\otimes (ts_2)^{-1}tg_1\right),$$
 which equals the first term.

 Let $\psi \in C^3(\C)^G$, we assert that $BA\psi = \psi.$ Recall that if $$f_3\otimes f_2\otimes f_1\in{}^{}_{u_4}\C[G]_{u_3}^{t_3} \otimes{}^{}_{u_3}\C[G]_{u_2}^{t_2}\otimes{}^{}_{u_2}\C[G]_{u_1}^{t_1},$$ then $$(A\psi)_{u_4,u_3,u_2,u_1}(f_3\otimes f_2\otimes f_1) = \psi (f_3\otimes t_3f_2\otimes t_3t_2f_1).$$
 Let $g_3\otimes g_2\otimes g_1 \in {}_{s_4u_4}\C_{s_3u_3}\otimes{}_{s_3u_3}\C_{s_2u_2}\otimes{}_{s_2u_2}\C_{s_1u_1}.$
 Then
 $$
 \begin{array}{lr}
      BA\psi(g_3\otimes g_2\otimes g_1)  =  s_4 A\psi(s_4^{-1}g_3\otimes s_3^{-1}g_2\otimes s_2^{-1}g_1)
       \end{array} $$

       where
$$s_4^{-1}g_3\otimes s_3^{-1}g_2\otimes s_2^{-1}g_1 \in{}_{u_4}\C[G]_{u_3}^{s_4^{-1}s_3} \otimes{}_{u_3}\C[G]_{u_2}^{s_3^{-1}s_2}\otimes{}_{u_2}\C[G]_{u_1}^{s_2^{-1}s_1}.$$
Hence
$$\begin{array}{llll}
BA\psi(g_3\otimes g_2\otimes g_1) = \\ s_4 \psi(s_4^{-1}g_3\otimes (s_4^{-1}s_3)s_3^{-1})g_2\otimes (s_4^{-1}s_3s_3^{-1}s_2)s_2^{-1}g_1)=\\
 s_4 \psi(s_4^{-1}g_3\otimes s_4^{-1}g_2\otimes s_4^{-1}g_1).
\end{array}
$$
Since $\psi$ is invariant, the later equals $\psi(g_3\otimes g_2\otimes g_1).$

Let $\varphi \in C_{\{1\}}(\C_T[G])$, next we will show $AB\varphi = \varphi.$ Consider
$$g_3\otimes g_2\otimes g_1\in{}_{t_4u_4}\C_{t_3u_3}\otimes {}_{t_3u_3}\C_{t_2u_2}\otimes{}_{t_2u_2}\C_{t_1u_1}.
$$
We have $$B\varphi (g_3\otimes g_2\otimes g_1)= t_4 \varphi(t_4^{-1}g_3\otimes t_3^{-1}g_2\otimes t_2^{-1} g_1).$$
Let $$f_3\otimes f_2\otimes f_1 \in {}^{}_{u_4}\C[G]^{s_3}_{u_3}\otimes {}^{}_{u_3}\C[G]^{s_2}_{u_2}\otimes {}^{}_{u_2}\C[G]_{u_1}^{s_1}.$$
Then
$$AB\varphi(f_3\otimes f_2\otimes f_1)= (B\varphi)(f_3\otimes s_3f_2\otimes s_3s_2f_1)$$ where
$$f_3\otimes s_3f_2\otimes s_3s_2f_1\in {}_{u_4}\C_{s_3u_3}\otimes{}_{s_3u_3}\C_{s_3s_2u_2}\otimes{}_{s_3s_2u_2}\C_{s_3s_2s_1u_1}$$
Hence
$$ AB\varphi(f_3\otimes f_2\otimes f_1)= \varphi (f_3\otimes s_3^{-1}s_3 f_2\otimes (s_3s_2)^{-1} s_3s_2 f_1) = \varphi(f_3\otimes f_2\otimes f_1).$$
\qed
\Edm

\subsubsection{\sf General case}
Our next aim is to show that the isomorphism of Theorem \ref{cohomology free} remains valid when the action of the group is not necessarily free. The following result has been proved in \cite{bw, herscovich}, see also \cite{hsjpaa}.
\Bpo
Let $\C$ and $\D$ be $k$-categories and let $F:\C\to \D$ be an equivalence of $k$-categories. There is an induced map
$$C^\bullet F: C^\bullet(\D)\to \C^\bullet(\C)$$
which is a quasi-isomorphism.
\Epo
 \Brm\label{future}
In the following the explicit definition of $C^\bullet( F )$ will be useful, it is as follows.
Let $$\varphi =\left(\varphi_{y_{n+1},\dots ,y_1}\right) \in C^n (\D) $$ where
    $$\varphi_{y_{n+1},\dots, y_1}\ :\ {}_{y_{n+1}}\D_{y_n}\otimes \cdots \otimes {}_{y_2}\D_{y_1} \longrightarrow {}_{y_{n+1}}\D_{y_1}$$ is a $k$-morphism.
The component $(x_{n+1},\cdots,x_1)$ of $(C^\bullet F)(\varphi)$ is given as follows. Let
   $$f_{n+1}\otimes\cdots \otimes f_1 \in {}_{x_{n+1}}\C_{x_n}\otimes \cdots \otimes {}_{x_2}\C_{x_1}.$$
   Then
   $$[(C^\bullet F)(\varphi)]_{x_{n+1}, \dots, x_1}(f_{n+1}\otimes\cdots \otimes f_1)= $$
   $$(_{x_{n+1}}F_{x_1})^{-1} \left(\varphi_{_{F(x_{n+1}), \dots, F(x_1)}}\left(F(f_{n+1})\otimes \cdots \otimes F(f_1)\right)\right)$$
   where $$_{x_{n+1}}F_{x_1}: {}_{x_{n+1}}\C_{x_1}\to {}_{F(x_{n+1})}\D_{F(x_1)}$$ is the $k$-isomophism provided by the equivalence $F.$
   \Erm

   Observe that in \cite{herscovich} the above Proposition is obtained in a  more general setting, that is for  a $\D$-bimodule of coefficients $\N$. In our case $\N=\D.$ The restricted $\C$-bimodule of coefficients is denoted $F\N$ in \cite{herscovich}, observe that $F\D$ is isomorphic to $\C$ via $F.$ This later isomorphism explains that in our setting ${}_{x_{n+1}}F_{x_1}^{-1}$ is required in the above formula while in \cite{herscovich} it is not needed since the bimodule of coefficients there is $F\D$, not $\C.$

\Bte
Let $\C$ and $\D$ be  $G$-$k$-categories and let $F:\C\to \D$ be a $G$-$k$- equivalence of categories. Then $F$ induces an isomorphism of $G$-$k$-algebras $$HH^\bullet (\D)\to HH^\bullet(\C).$$
\Ete
\Bdm
The explicit description of $C^\bullet F$ given above enables to check without difficulty that it is multiplicative with respect to the cup product.
Moreover, $C^\bullet F$ commutes with the actions of $G$ on $C^\bullet\C$ and $C^\bullet \D,$ that is $C^\bullet F$ is a $kG$-morphism. Therefore the induced map in cohomology is an isomorphism of $G$-$k$-algebras.
\qed
\Edm
We recall that if $\C$ is a $G$-$k$-category, then $M_G(\C)$ is a $G$-$k$-category  where the action of $G$ on the objects of $M_G(\C)$ is free, see Definition \ref{MGC}. Moreover there is a  a $G$-$k$-functor $L:M_G(\C) \to \C$ which is an equivalence of categories.
 \Bte\label{theorem comparison cohomology general case}
 Let $\C$ be a $G$-$k$-category. Let $\C[G]$ be the graded skew category, and let $\{1\}$ be the trivial conjugacy class of $G$.
 There is an isomorphism of $k$-algebras
 $$HH^*_{\{1\}} (\C[G]) \simeq H^*(C^\bullet(\C)^G).$$
 If the invariants  functor  is exact, we have an isomorphism of $k$-algebras

  $$HH^*_{\{1\}} (\C[G]) \simeq HH^*(\C)^G.$$
  \Ete
\Bdm
Let $$L[G]: M_G(\C)[G] \to \C[G]$$ be the homogeneous equivalence of $G$-graded $k$-categories obtained from the $G$-$k$-equivalence of categories
$L:M_G(\C)\to \C$ of Theorem \ref{L}.

We observe that if $\B$ and $\D$ are $G$-graded categories and $K:\B\to\D$ is a homogeneous equivalence, then the quasi-isomorphism
$$C^\bullet(K): C^\bullet(\B) \to C^\bullet(\D)$$ described in Remark \ref{future} preserves the decomposition along the conjugacy classes of
$G.$ Hence $HH^*_{\{1\}}(\B)$ and $HH^*_{\{1\}}(\D)$ are isomorphic $k$-algebras. \qed

\Edm

\subsection{\sf Galois coverings}\label{section Galois cohomology}

In this Subsection we will translate the results we have obtained for the cohomology of skew categories to a Galois coverings $\C\to \C/G$. Then we will provide a canonical monomorphism from the invariants of the cohomology of $\C/G$ to the cohomology of $\C$. This corresponds to the  monomorphism of \cite{mmm}. In low degrees it is described in \cite{ghs}.

The proof of the following result is along the same lines than the proof of Theorem \ref{theorem cohomology Galois covering}.

\Bte
Let $\C \to \C/G$ be a Galois covering.
$$HH^*_{\{1\}}(\C/G) = H_*(C_\bullet(\C)^G).$$
If the invariants functors is exact, then
$$HH^*_{\{1\}}(\C/G) = HH^*(\C)^G.$$
\Ete

\Bcr
Let $\C \to \C/G$ be a Galois covering. If the invariants functor is exact, there is a canonical injective morphism
$$HH^*(\C)^G \hookrightarrow HH^*(\C/G)$$
which splits canonically.
\Ecr

\Bdm
The cohomology of $\C/G$ has a direct sum decomposition along the conjugacy classes of $G$. The direct summand corresponding to the trivial conjugacy class is isomorphic to the invariants of the cohomology of $G$. \qed
\Edm

\section{\sf Hochschild cohomology of skew group algebras}\label{section skew group algebras}

In this section we will specialise Theorem \ref{theorem comparison cohomology general case} for $k$-algebras. Note that the proof of Theorem \ref{theorem comparison cohomology general case}  requires the resolving category. We do not know a proof of Theorem \ref{theo for skew group algebras} without using a resolving object which  makes the action of the group free on a set. See Remark \ref{final remark}.

Let $\Lambda$ be a $k$-algebra, and let $G$ be a group acting by algebra automorphisms of $\Lambda$.  Let $\Lambda[G]$ be the usual skew group algebra recalled in Remark \ref{definition skew group algebra}.
The Hochschild cohomology of a $k$-algebra $\Lambda$  is denoted $HH^*(\Lambda)$.

\begin{theo}\label{theo for skew group algebras}
Let $G$ be a finite group whose order is invertible in $k$. Let $\Lambda$ be a $k$-algebra with an action of $G$ by algebra automorphisms. There is an isomorphism of algebras
$$HH_{\{1\}}^*(\Lambda[G]) \simeq HH^*(\Lambda)^G.$$
\end{theo}

\Bdm
Let $\Lambda_1$ be the single object $G$-$k$-category of $\Lambda$ considered at Example \ref{single object}. As noticed, the action of $G$ is not free on $\Lambda_1$ unless $G$ is trivial. By Theorem \ref{theorem comparison cohomology general case}  we have an isomorphism of $k$-algebras
 \begin{equation}\label{un}
   HH_{\{1\}}^*(\Lambda_1[G]) \simeq HH^*(\Lambda_1)^G.
 \end{equation}

  Of course $a(\Lambda_1)=\Lambda$. Moreover  we have that $HH^*(\C)=HH^*(a(\C))$, see for instance \cite{cr}. Hence the right hand side of (\ref{un}) is isomorphic to $HH^*(\Lambda)^G$.

 On the other hand, as quoted in the Introduction, if $G$ is finite and if $\C$ is a $G$-$k$-category with a finite number of objects, then
 $$a(\C[G])= a(\C)[G].$$
 Thus the left hand side of (\ref{un}) is
 $$HH_{\{1\}}^*(\Lambda_1[G])=HH_{\{1\}}^*(a(\Lambda_1[G]))= HH_{\{1\}}^*(a(\Lambda_1)[G])= HH_{\{1\}}^*(\Lambda[G]).$$\qed

\Edm

\begin{rema}\label{final remark}
 The proof above of Theorem \ref{theo for skew group algebras} relies on Theorem \ref{theorem comparison cohomology general case}. Consequently, it is interesting to track the proof of Theorem \ref{theorem comparison cohomology general case} specifying it to an algebra. First, we have considered the matrix algebra $M_G(\Lambda)$ described in Proposition \ref{resolving algebra of an algebra}. Then the categorical proof translates into  decomposing the cochains of $M_G(\Lambda)$ through the set $E$ of diagonal idempotents of Remark \ref{free idempotents}. The  freeness of the action on $E$ enables to show that the module of invariants of the complex of cochains of $M_G(\Lambda)$ is isomorphic to the homogeneous cochains of the conjugacy class $1$ of $M_G(\Lambda)[G]$. The final step consists in showing that the Hochschild cohomology of $\Lambda[G]$, as a $kG$-module, remains the same when considering the  algebra $M_G(\Lambda)[G]$.
\end{rema}

\footnotesize
\noindent C.C.:\\
IMAG, Univ Montpellier, CNRS, Montpellier, France\\
Institut Montpelli\'{e}rain Alexander Grothendieck  \\{\tt Claude.Cibils@umontpellier.fr}

\medskip

\noindent E.N.M.:\\
Departamento de Matem\'atica, IME-USP,\\
Rua do Mat\~ao 1010, cx pt 20570, S\~ao Paulo, Brasil.\\
{\tt enmarcos@ime.usp.br}


\begin{thebibliography}{99}
\bibitem{abm} Assem, I.; Bustamante, J. C.; Le Meur, P.; Coverings of Laura Algebras: the Standard Case.
J. Algebra 323 (2010) 83--120.



\bibitem{bw} Baues H.J.; Wirsching G.; Cohomology of small categories. J. Pure Appl. Algebra 38 (1985), 187--211.

\bibitem{bg} Bongartz; K., Gabriel, P.; Covering spaces in representation-theory. Invent. Math. 65 (1982), 331--378.

\bibitem{ce} Cartan, H.; Eilenberg S.; Homological Algebra, Princeton, New Jersey, Princeton University press, 1958.

\bibitem{com} Cohen, M.; Montgomery, S.; Group-graded rings, smash products, and group actions. Trans. Am. Math. Soc. 282 (1984), 237--258.

%\bibitem{cc} Cibils, C.; Cyclic and Hochschild homology of 2-nilpotent algebras. K-Theory 4 (1990), 131--141.

\bibitem{cm} Cibils, C.; Marcos, E. N.; Skew category, Galois covering and smash product of a k-category.
Proc. Amer. Math. Soc. 134 (2006), 39--50.



\bibitem{cr} Cibils, C.; Redondo, M. J.; Cartan-Leray spectral sequence for Galois coverings of linear categories. J. Algebra 284 (2005), 310--325.

\bibitem{cs} Cibils, C.; Solotar, A.; The fundamental group of a Hopf linear category. J. Algebra 462 (2016), 137--162.

\bibitem{cornick} Cornick, J.; Homological techniques for strongly graded rings: a survey. Geometry and cohomology in group theory (Durham, 1994),
88--107, London Math. Soc. Lecture Note Ser., 252, Cambridge Univ. Press, Cambridge, 1998.

\bibitem{coell} Corti\~{n}as, G.; Ellis, E.; Isomorphism conjectures with proper coefficients.
J. Pure Appl. Algebra 218 (2014), 1224--1263.

\bibitem{farinati} Farinati M.; Hochschild duality, localization, and smash products. J. Algebra 284 (2005), 415--434.

\bibitem{gabriel} Gabriel, P.; The universal cover of a representation-finite algebra. Representations of algebras (Puebla, 1980), pp. 68--105, Lecture Notes in Math., 903, Springer, Berlin-New York, 1981.

\bibitem{gp} Garc\'ia, O. C.; de la Pe\~na, J.A.; Lattices with a finite Whitman cover. Algebra Univers. 16 (1983), 186--194.


\bibitem{gk} Ginzburg, V.; Kaledin, D.: Poisson deformations of symplectic quotient singularities, Adv. Math. 186 (2004), 1--57.

\bibitem{ghs} Green, E.L,; Hunton J. ; Snashall, N.; Coverings, the graded center and
Hochschild cohomology. J. Pure Appl. Algebra 212 (2008), 2691--2706.

\bibitem{herscovich} Herscovich, E.; La homolog\'\i a de Hochschild-Mitchell de categor\'\i as lineales y su parecido con la homolog\'\i a de Hochschild de \'algebras. Tesis de licenciatura en Ciencias Matem\'aticas,
Facultad de Ciencias Exactas y Naturales, Universidad de Buenos Aires (2005).\nolinebreak
\footnotesize
\begin{verbatim}
https://www-fourier.ujf-grenoble.fr/~eherscov/Articles/MasMath.pdf
\end{verbatim}

\normalsize
\bibitem{hsjpaa} Herscovich, E.; Solotar, A.; Hochschild-Mitchell cohomology and Galois extensions. J. Pure Appl. Algebra 209 (2007),  37--55.


\bibitem{hsja}Herscovich, E.; Solotar, A.;
Derived invariance of Hochschild-Mitchell (co)homology and one-point extensions.
J. Algebra 315 (2007),  852--873.

\bibitem{kp} Kasjan S.; De la Pe\~na, J. A.; Galois coverings and the problem of axiomatization of the representation type of algebras. Extracta Math. 20 (2005), 137--150.

\bibitem{keller} Keller, B.;  On the cyclic homology of exact categories. J. Pure Appl. Algebra 136 (1999), 1--56.



\bibitem{loday}Loday J. L.; Cyclic homology. 2nd ed. Grundlehren der Mathematischen Wissenschaften. 301. Berlin: Springer. 513 p. (1998).

\bibitem{lorenz} Lorenz, M.;  On Galois descent for Hochschild and cyclic homology. Comment. Math. Helv. 69 (1994), 474–482.

\bibitem{mccarthy} McCarthy, R.; The cyclic homology of an exact category. J. Pure Appl. Algebra 93 (1994), 251--296.


\bibitem{mmm} Marcos, E.; Martinez-Villa, R.; Martins, M. I.; Hochschild cohomology of skew group rings and invariants. Cent. Eur. J. Math. 2 (2004), 177--190.

    \bibitem{mitchell} Mitchell, B.; Rings with several objects. Advances in Math. 8 (1972), 1--161.

 \bibitem{pr}   Pirashvili, T.; Redondo, M.J.
Cohomology of the Grothendieck construction. Manuscr. Math. 120 (2006), 151--162.

\bibitem{nsw} Naile, D.;  Shroff, P.; Witherspoon,S.; Hochschild cohomology of group extensions of quantum symmetric algebras. Proc. Amer. Math. Soc. 139 (2011), 1553--1567.

\bibitem{ri} Riedtmann, C.; Algebren, Darstellungsk\"{o}cher, \"{U}berlagerungen und zur\"{u}ck. Comment. Math. Helv. 55 (1990), 199--224.

\bibitem{sw} Shepler A.; Witherspoon S.; Finite groups acting linearly: Hochschild cohomology and the cup product. Adv. Math. 226 (2011), 2884--2910.

\bibitem{stefan} Stefan D.; Hochschild cohomology on Hopf Galois extensions. J. Pure Appl. Algebra 103 (1995), 221--233.

\bibitem{stefan1996} Stefan D.; Decomposition of Hochschild homology of Hopf Galois extensions. Commun. Algebra 24 (1996), 1695--1706.

\bibitem{weibel} Weibel, C.; An introduction to homological algebra.
Cambridge Studies in Advanced Mathematics. 38. Cambridge: Cambridge University Press. 450 p. (1994).





\end{thebibliography}
\end{document}